\documentclass[review]{elsarticle}

\usepackage{lineno,hyperref}
\modulolinenumbers[5]

\journal{Systems \& Control Letters }

\usepackage[utf8]{inputenc}
\usepackage{amsmath}
\usepackage{amsfonts}
\usepackage{amssymb}
\usepackage{theorem,pifont}
\usepackage{graphicx}
\usepackage{color}

\newtheorem{theorem}{Theorem}
\newtheorem{corollary}{Corollary}

\newtheorem{definition}{Definition}
\newtheorem{example}{Example}
\newtheorem{remark}{Remark}
\newtheorem{lemma}{Lemma}
\newtheorem{proposition}{Proposition}

\newcommand{\mR}{\mathbb R}
\newcommand{\mRp}{\mathbb R_{\geq 0}}

\begin{document}

\begin{frontmatter}

\title{Growth conditions for global exponential stability and exp-ISS of time-delay systems under point-wise dissipation}

\author[antoine]{Antoine Chaillet\corref{mycorrespondingauthor}}
\cortext[mycorrespondingauthor]{Corresponding author}
\ead{antoine.chaillet@centralesupelec.fr}

\author[iasson]{Iasson Karafyllis}

\author[pp]{Pierdomenico Pepe}

\author[yuan]{Yuan Wang}

\address[iasson]{Mathematics Department, National Technical University of Athens, Athens 15780, Greece.}
\address[antoine]{Universit\'e Paris-Saclay, CNRS, CentraleSup\'elec, Laboratoire des signaux et syst\`emes, 91190, Gif-sur-Yvette, France.}
\address[pp]{Information Engineering, Computer Science, and Mathematics Department, University of L'Aquila, 67100, L'Aquila, Italy.}
\address[yuan]{Mathematical Sciences Department, Florida Atlantic University, Boca Raton, FL 33431, U.S.A.}

\begin{abstract}
For time-delay systems, it is known that global asymptotic stability is guaranteed by the existence of a Lyapunov-Krasovskii functional that dissipates in a point-wise manner along solutions, namely whose dissipation rate involves only the current value of the solution's norm. So far, the extension of this result to global exponential stability (GES) holds only for systems ruled by a globally Lipschitz vector field and remains largely open for the input-to-state stability (ISS) property. In this paper, we rely on the notion of exponential ISS to extend the class of systems for which GES or ISS can be concluded from a point-wise dissipation. Our results in turn show that these properties still hold in the presence of a sufficiently small additional term involving the whole state history norm. We provide explicit estimates of the tolerable magnitude of this extra term and show through an example how it can be used to assess robustness with respect to modeling uncertainties.
\end{abstract}


\end{frontmatter}


\section{Introduction}

The input-to-state stability (ISS) framework has become a central and classical tool to study stability and robustness of nonlinear systems. Originally developped for systems described by ordinary differential equations \cite{SON1,cetraro}, it has progressively been extended to infinite-dimensional systems \cite{mironchenko2020input,KAKR19book}.

A particular class of infinite-dimensional systems is that of time-delay systems. Due to the peculiarities of this subclass, it has been the subject of specific ISS developments initiated with the works \cite{Teel98,Pepe:2006ju}. ISS for time-delay systems has now become a mature topic, but some fundamental questions remain open.

In particular, a Lyapounov-Krasovskii functional (LKF) characterization of ISS was provided in \cite{Karafyllis:2008hc,Kankanamalage:2017ug}. This characterization requests that the dissipation of the LKF along the system's solutions is expressed in terms of the LKF itself (\emph{LKF-wise dissipation}). This requirement turns out to be rather unhandy in practice. More crucially, it is not in line with the LKF characterization of global asymptotic stability, in which the LKF is allowed to dissipate merely in terms of the current value of the solution's norm (\emph{point-wise dissipation)} \cite{KRA59}. To date, it is not known whether ISS can be ensured through a point-wise dissipation \cite{CHPEMACH17}, although this question has recently received a positive answer for the weaker notion of integral ISS \cite{CHGOPE22}.

Even in the absence of exogenous inputs, point-wise dissipation is not yet fully understood. In particular, the possibility to ensure global exponential stability with a point-wise dissipation still constitutes an open question. A partial answer was given in \cite{chaillet2019relaxed}, but only for systems ruled by a globally Lipschitz vector field. A positive answer to this question would complement the existing arsenal to establish GES, including the Razumikhin approach \cite{wang2005exponential} and Halanay's inequality \cite[p. 378]{Halanay:1966wl}.

In this paper, we significantly enlarge the class of systems for which ISS or GES can be established under a point-wise dissipation. To that aim, we focus on the exp-ISS property, which is a particular case of ISS in which the influence of the initial state is requested to decay exponentially. In particular, exp-ISS ensures GES for the corresponding input-free system. We provide growth conditions under which a point-wise dissipation is enough to conclude exp-ISS. These conditions may take two forms: either left or right, depending on whether the increase rate or the decay rate is restricted. Both these conditions turn out to be automatically satisfied when the vector field is globally Lipschitz, but we show through examples that the proposed growth conditions encompass a much wider class of systems, thus significantly generalizing the results in \cite{chaillet2019relaxed}.

While making the analysis simpler, the use of a point-wise dissipation complicates the analysis of robustness to modeling errors or parameter uncertainties, as compared to an LKF-wise dissipation or a history-wise one (in which the LKF is requested to dissipate in terms of the whole state history norm). To address this issue, we allow for an additional quadratic positive term in the LKF's derivative. This term involves the full state history. Since the dissipation is only in terms of the current solution's norm, the derivative of the considered LKF is thus no longer guaranteed to be non-positive, even in the absence of an input. Despite this severe issue, we show that exp-ISS does hold provided that this extra term is sufficiently small. Our proofs being constructive, we actually provide explicit estimates of the strength of this additional quadratic term. We show through an example how this result can indeed be useful for robustness analysis with respect to modeling errors.

On our way to establish these results, we provide several technical lemmas that may be of interest on their own. In particular, we propose a novel sufficient condition for a property known as \emph{robust forward completeness} or \emph{bounded reachability property}, which plays a central role in the stability analysis of time-delay systems \cite{Mironchenko:2017wb, KAPECHWA22-hal}. We also show that exp-ISS can be reformulated as two specific inequalities on the system's solutions.

\vspace{3mm}
\noindent\textbf{Notation.} Given $n\in\mathbb N_{\geq 1}$ and $\Delta\geq 0$, $\mathcal X^n$  denotes the set of all continuous functions from $[-\Delta,0]$ to $\mR^n$, whereas $\mathcal W^n$ denotes the Sobolev space of absolutely continuous functions mapping $[-\Delta,0]$ into $\mathbb R^{n}$ with essentially bounded derivative, and $C^1([-\Delta,0];\mR^n)$ denotes the set of all continuously differentiable functions from $[-\Delta,0]$ to $\mR^n$. 
Given a continuous signal $x:[-\Delta,T)\to \mR^n$ with $T\in(0,+\infty]$ and any $t\in [0,T)$, $x_t\in\mathcal X^n$ denotes the history function defined as $x_t(\tau):=x(t+\tau)$ for all $\tau\in[-\Delta,0]$. Given a non-empty interval $I\subset\mR$ and an essentially bounded Lebesgue measurable signal $u:I\to\mR^m$, $\Vert u\Vert:=\textrm{ess\,sup}_{t\in I}|u(t)|$, where $|\cdot|$ denotes the Euclidean norm. Given $m\in\mathbb N_{\geq 1}$, $\mathcal U^m$ denotes the set of all signals $u:\mRp\to\mR^m$ that are Lebesgue measurable and locally essentially bounded. Given $u\in\mathcal U^m$ and $t_1\geq t_2\geq 0$, $u_{[t_1,t_2]}$ denotes the restriction of $u$ to $[t_1,t_2]$, namely $u:[t_1,t_2]\to\mR^m$ is defined as $u_{[t_1,t_2]}(t):=u(t)$ for all $t\in[t_1,t_2]$. Given a continuously differentiable function $V:\mR^n\to\mR$, $\nabla V$ denotes its gradient. Given $a\in\mRp$, $\lceil a\rceil$ denotes the smallest $q\in\mathbb N$ for which $q\geq a$. A function $\alpha:\mRp\to\mRp$ is said to be of class $\mathcal N$ if it is continuous, non-decreasing and zero at zero. It is said of class $\mathcal K$ if, in addition, it is increasing. It is said to be of class $\mathcal K_\infty$ if it is of class $\mathcal K$ and satisfies $\lim_{s\to+\infty}\alpha(s)=+\infty$. Given a functional $V:\mathcal {X}^n\to \mathbb R_{\ge 0}$, its Driver derivative $D^+V:\mathcal X^n\times \mR^n\to[-\infty,+\infty]$ is defined, for all $\phi \in \mathcal X^n$ and all $w\in \mathbb R^n$, as
$D^+ V(\phi,w):=
\limsup_{h\to 0^+} \frac {V(\phi^{}_{h,w})-V(\phi)}{h}$ where, for each $h\in [0,\Delta)$, $\phi^{}_{h,w} \in \mathcal {X}^n$ is given by
\begin{equation*}
\phi^{}_{h,w} (s) :=  \left \{ \begin{array}{cl}  \phi(s+h), &\textrm{ if } s\in [-\Delta,-h),    \\
\phi(0)+(h+s)w, &\textrm{ if } s\in [-h,0].
\end{array} \right .
\end{equation*}

\section{Context}

\subsection{Global exponential stability}

We start by considering input-free delay systems, namely:
\begin{align}\label{eq-0}
\dot x(t)=f_0(x_t),
\end{align}
where $x_t\in\mathcal X^n$ and $f_0:\mathcal X^n\to\mR^n$ is a vector field which is Lipschitz on bounded sets and satisfies $f_0(0)=0$. For such class of systems, we recall the definition of global exponential stability.

\begin{definition}[GES]
The origin of \eqref{eq-0} is said to be \emph{globally exponentially stable (GES)} if there exist $k,\eta>0$ such that, for all $x_0\in\mathcal X^n$, the corresponding solution of \eqref{eq-0} satisfies $|x(t)|\leq k\|x_0\|e^{-\eta t}$ for all $t\geq 0$.
\end{definition}

The GES property therefore ensures that all solutions converge at the exponential rate $\eta$ and that their transient overshoot is bounded by a linear function of the initial state norm. A powerful tool to establish GES is through a Lyapunov-Krasovskii functional (LKF). In particular, we recall the following from \cite{KRA59,Pepe:2013it}.

\begin{theorem}[Existing LKF characterization for GES]\label{theo-GES-char} The following statements are equivalent:
\begin{itemize}
\item[$i)$] the origin of \eqref{eq-0} is GES
\item[$ii)$] there exist a functional $V:\mathcal {X}^n\to \mathbb R_{\ge 0}$, Lipschitz on bounded sets, and $\underline a,\overline a, a, \rho>0$ such that, for all $\phi \in \mathcal {X}^n$, 
\begin{align*}
\underline a|\phi(0)|^\rho\le V(\phi)\le \overline a\|\phi\|^\rho \\
D^+V(\phi,f_0(\phi))\le -aV(\phi)
\end{align*}
\item[$iii)$] there exist a functional $V:\mathcal {X}^n\to \mathbb R_{\ge 0}$, Lipschitz on bounded sets, and $\underline a,\overline a, a>0$ such that, for all $\phi \in \mathcal {X}^n$, 
\begin{align*}
\underline a\|\phi\|\le V(\phi)\le \overline a\|\phi\| \\
D^+V(\phi,f_0(\phi))\le -a\|\phi\|.
\end{align*}
\end{itemize}
\end{theorem}

In this statement, $D^+V(\phi,f(\phi,v))$ denotes Driver's derivative of $V$ along the system's solutions, as recalled in the Notation paragraph. In particular, $D^+V(x_t,f(x_t,u(t))$ coincides almost everywhere with the upper-right Dini derivative of $t\mapsto V(x_t)$ on its maximal interval of existence: see \cite{driver1962existence} and \cite[Theorem 2]{pepe2007liapunov}. 

Beyond the fact that $\rho$ does not need to be equal to 1 in Item $ii)$, there are two key differences between Items $ii)$ and $iii)$. The first one lies in the fact that the LKF $V$ in Item $iii)$ is \emph{coercive}, in the sense that it can only vanish if the whole state history is identically zero. On the contrary, the LKF in Item $ii)$ is lower-bounded only in terms of $|\phi(0)|$, meaning that it may vanish even if $\phi\neq 0$. The possibility offered by Item $ii)$ to consider non-coercive LKFs turns out to be very convenient in practice, whereas the fact that $V$ is coercive in Item $iii)$ often proves useful to conduct further robustness analysis. 

The second difference is in the way the LKF dissipates along solutions: in Item $ii)$, this dissipation is in terms of the LKF itself (\emph{LKF-wise dissipation}), whereas, in Item $iii)$, the dissipation is requested to involve the whole history norm (\emph{history-wise dissipation}). When the LKF is coercive, there is no qualitative difference between these two types of dissipation, but LKF-wise dissipation does not necessarily guarantee a history-wise dissipation when $V$ is not coercive. 

To sum up, while Item $ii)$ is more convenient to establish GES, Item $iii)$ constitutes a powerful converse result once GES is known to hold.

\subsection{Exponential ISS}

We may consider the impact of an exogenous input on the system, which then takes the form
\begin{align}\label{eq-01}
\dot x(t)=f(x_t,u(t)),
\end{align}
where $x_t\in\mathcal X^n$, $u\in\mathcal U^m$, and $f:\mathcal X^n\times\mR^m\to\mR^n$ is a functional which is Lipschitz on bounded sets and satisfies $f(0,0)=0$. A powerful way to assess robustness with respect to the input signal $u$ is through the input-to-state stability (ISS) framework. The ISS property requests in particular that the steady-state error of solutions is small if the applied input is of sufficiently small amplitude \cite{cetraro,Pepe:2006ju}. In this paper, we focus on the following particular declination of ISS.

\begin{definition}[Exp-ISS, linear gain]
The system \eqref{eq-01} is said to be \emph{exponentially input-to-state stable (exp-ISS)} if there exist $k,\eta>0$ and $\mu\in\mathcal N$ such that, for all $x_0\in\mathcal X^n$ and all $u\in\mathcal U^m$, the corresponding solution of \eqref{eq-01} satisfies
\begin{align*}
|x(t)|\leq k\|x_0\|e^{-\eta t}+\mu(\|u_{[0,t]}\|),\quad \forall t\geq 0.
\end{align*}
It is said to be \emph{exp-ISS with linear gain} if, in addition, there exists $\mu_0\geq 0$ such that $\mu(s)=\mu_0s$ for all $s\geq 0$.
\end{definition}

This property has already been used in the literature of infinite-dimensional systems. In particular, it was studied in \cite{mironchenko2019monotonicity} in the context of parabolic partial differential equations. Before that, it was used for general infinite-dimensional systems under the name \emph{eISS} \cite{Dashkovskiy:2012hq}. It should not be confused with the exponential ISS notion employed in \cite{PRAWAN}, which rather considers ISS with respect to a filtered version of the input.

The above exp-ISS property clearly ensures 0-GES, meaning that the origin of the input-free system $\dot x(t)=f(x_t,0)$ is GES. But it also ensures that, in response to a bounded input $u$, any solution is attracted by a $\mu(\|u\|)-$neighborhood of the origin. In particular, solutions are bounded in response to any bounded input, the steady-state error is small for sufficiently small inputs, and solutions converge to the origin in response to any vanishing input. In the particular case of exp-ISS with linear gain, the steady-state error is at most proportional to the amplitude of the applied input. This feature mimics what happens for linear systems and turns out particularly useful when invoking small-gain results for the stability analysis of interconnected systems \cite{Karafyllis:2007kz}.

Exp-ISS can be established by existing LKF tools. In particular, the following result can be proved using classical manipulations.

\begin{proposition}[LKF-wise dissipation for exp-ISS]\label{prop-existing-LKF-expISS}
Assume that there exist a functional $V:\mathcal X^n\to\mRp$ which is Lipschitz on bounded sets, constants $\underline a,\overline a,a,\rho>0$ and a function $\gamma\in\mathcal N$ such that, for all $\phi\in\mathcal X^n$ and all $v\in\mR^m$,
\begin{align}
\underline a |\phi(0)|^\rho&\leq V(\phi)\leq \overline a \|\phi\|^\rho\\
D^+V(\phi,f(\phi,v))&\leq -aV(\phi)+\gamma(|v|).
\end{align}
Then \eqref{eq-01} is exp-ISS. Moreover, if there exists $\gamma_0\geq 0$ such that $\gamma(s)=\gamma_0s^\rho$ for all $s\geq 0$, then \eqref{eq-01} is exp-ISS with linear gain.
\end{proposition}

Just like Theorem \ref{theo-GES-char}, this result requests a LKF-wise dissipation. Interestingly, Proposition \ref{prop-existing-LKF-expISS} can actually be strengthened by allowing a positive square term of the state history norm in the expression of the LKF's derivative, as stated next.

\begin{theorem}[Alternative LKF-wise dissipation for exp-ISS]\label{theo-LKF-expISS-nogrowth}
Assume that there exist a functional $V:\mathcal X^n\to\mRp$ which is Lipschitz on bounded sets, constants $\underline a,\overline a,a,\rho>0$, $c\geq 0$, and a function $\gamma\in\mathcal N$ such that, for all $\phi\in\mathcal X^n$ and all $v\in\mR^m$,
\begin{align}
\underline a |\phi(0)|^\rho &\leq V(\phi)\leq \overline a \|\phi\|^\rho\label{eq-57}\\
D^+V(\phi,f(\phi,v))&\leq -aV(\phi)+c\|\phi\|^\rho+\gamma(|v|).\label{eq-58}
\end{align}
Then, provided that
\begin{align}\label{eq-56}
c< \underline a ae^{-a\Delta},
\end{align}
the system \eqref{eq-01} is exp-ISS . If, in addition, there exists $\gamma_0\geq 0$ such that $\gamma(s)=\gamma_0s^\rho$ for all $s\geq 0$, then \eqref{eq-01} is exp-ISS with linear gain.
\end{theorem}

This result, proved in Section \ref{sec-proof-theo-LKF-expISS-nogrowth}, suggests some robustness of the exp-ISS property: as compared to Proposition \ref{prop-existing-LKF-expISS}, the additional positive term $c\|\phi\|^\rho$ does not compromise exp-ISS provided that $c$ is small enough to satisfy \eqref{eq-56}. This would come as no surprise if the considered LKF was coercive (meaning lower bounded by a term of the form $\underline a \|\phi\|^\rho$), which is however not requested here. Note that, even in the input-free case, this result seems to be novel and provides a useful sufficient condition for GES (simply consider $v=0$ in Theorem \ref{theo-LKF-expISS-nogrowth}).

\section{Growth conditions for exp-ISS}
 
While Theorems \ref{theo-GES-char} and \ref{theo-LKF-expISS-nogrowth} provide powerful conditions for GES and exp-ISS, the requirement of a LKF-wise (or history-wise) dissipation often significantly complicates the analysis. It also lacks homogeneity with the existing LKF condition for global asymptotic stability, in which a point-wise dissipation is enough \cite{KRA59}. Here, we investigate whether such a \emph{point-wise} dissipation would be sufficient to ensure GES or exp-ISS. 

Partial answers in this direction were given in \cite{CHPEMACH17}, where growth conditions were provided to ensure that a point-wise dissipation ensures ISS. For the GES property, it was shown in \cite{chaillet2019relaxed} that such a point-wise dissipation is indeed sufficient provided that the vector field is globally Lipschitz, which constitutes a conservative constraint in practice. To date, it is not known whether a point-wise dissipation can be used to establish GES or exp-ISS, but the results presented below enlarge the class of systems for which this holds true by imposing specific growth conditions. These conditions take two possible forms: either right or left, depending on whether we assume a sufficiently slow increase or decrease. We detail these results in two separate sections.

\subsection{Right growth condition}

Our first result ensures exp-ISS under a point-wise dissipation provided an upper bound on a specific function involving the vector field. In line with Theorem \ref{theo-LKF-expISS-nogrowth}, we additionally allow for the presence of a positive term in the state history norm, but the result is new also when this term is taken as zero. Its proof is reported in Section \ref{sec-proof-theo-expISS2}.

\begin{theorem}[Exp-ISS under an alternative right growth condition]\label{theo-expISS2}
Assume that there exist a functional $V:\mathcal X^n\to\mathbb R_{\geq 0}$ which is Lipschitz on bounded sets, $\overline a,a>0$, $c\geq 0$, and $\gamma\in\mathcal N$  such that, for all $\phi\in\mathcal X^n$ and all $v\in\mathbb R^m$,
\begin{align}
0\leq V(\phi)&\leq \overline a\|\phi\|^2\label{eq-2}\\
D^+V(\phi,f(\phi,v))&\leq -a|\phi(0)|^2+c\|\phi\|^2+\gamma(|v|).\label{eq-3}
\end{align}
Assume further that there exists a symmetric positive definite matrix $P\in\mR^{n\times n}$ and a constant $\sigma>0$ such that, for all $\phi\in\mathcal X^n$ and all $v\in\mR^m$,
\begin{align}
\phi(0)^\top Pf(\phi,v)&\leq \sigma(\|\phi\|^2+\gamma(|v|)).\label{eq-22}
\end{align}
Then there exists $\bar c>0$ such that, if $c\in[0,\bar c)$, the system \eqref{eq-01} is exp-ISS. If, in addition, there exists $\gamma_0\geq 0$ such that $\gamma(s)=\gamma_0s^2$ for all $s\geq 0$, then \eqref{eq-01} is exp-ISS with linear gain.
\end{theorem}


This result provides a right growth condition on the vector field $f$ under which a point-wise dissipation is enough to guarantee exp-ISS. It is worth stressing that condition \eqref{eq-22} holds automatically for any arbitrary positive definite matrix $P$ when $f$ is globally Lipschitz in the state, uniformly in the input, namely if there exists $\ell\geq 0$ such that 
\begin{align*}
|f(\phi,v)-f(\psi,v)|\leq \ell \|\phi-\psi\|, \quad \forall \phi,\psi\in\mathcal X^n,\ v\in\mR^n.
\end{align*}
In this case, \eqref{eq-22} turns out to be satisfied with any $\gamma\in\mathcal N$ satisfying $|f(0,v)|^2\leq \gamma(|v|)$ for all $v\in\mR^m$.
Nevertheless, as can be seen from Examples \ref{exa-1}-\ref{exa-3} below, the proposed growth condition is far from being restricted to such class of systems, thus significantly generalizing \cite[Corollary 1]{chaillet2019relaxed}. The requirement \eqref{eq-22} can be interpreted as a condition to guarantee that the current history norm $\|x_t\|$ is not too much greater than the current solution norm $|x(t)|$. More precisely, it imposes a restriction on the increase rate of the function $t\mapsto x(t)^\top Px(t)$.  

Rather surprisingly, the above result does not impose any lower bound on the considered LKF, not even in terms of $|\phi(0)|$: see \eqref{eq-2}. In particular, it allows for non-coercive LKFs, which is precious from an application viewpoint.

An immediate consequence of Theorem \ref{theo-expISS2} is the following sufficient condition for GES.

\begin{corollary}[GES under a right growth condition]\label{corol-GES}
Assume that there exist a functional $V:\mathcal X^n\to\mathbb R_{\geq 0}$ which is Lipschitz on bounded sets, $\overline a,a>0$ and $c\geq 0$ such that, for all $\phi\in\mathcal X^n$, \eqref{eq-2} holds and
\begin{align*}
D^+V(\phi,f_0(\phi))&\leq -a|\phi(0)|^2+c\|\phi\|^2.
\end{align*}
Assume further that there exists a symmetric positive definite matrix $P\in\mR^{n\times n}$ and a constant $\sigma>0$ such that, for all $\phi\in\mathcal X^n$,
\begin{align*}
\phi(0)^\top Pf_0(\phi)&\leq \sigma \|\phi\|^2.
\end{align*}
Then there exists $\bar c>0$ such that, if $c\in[0,\bar c)$, the origin of \eqref{eq-0} is GES.
\end{corollary}

\begin{remark}[Estimate of $\bar c$]\label{rmk1}
The proof  shows that $\bar c$ can be picked in Theorem \ref{theo-expISS2} and Corollary \ref{corol-GES} as\footnote{See \eqref{eq-29} and \eqref{eq-30}.}
\begin{align*}
\bar c=\min\left\{ \frac{p_m}{\sigma}e^{-2\Delta}\,,\,1\right\}\frac{a p_m}{2p_M}e^{-2\Delta},
\end{align*}
where $p_m,p_M>0$ denote the smallest and largest eigenvalues of $P$ respectively. This may however constitute a conservative estimate of the maximal positive quadratic term appearing in Driver's derivative of the considered LKF.
\end{remark}

The proof of Theorem \ref{theo-expISS2} consists in explicitly constructing a coercive LKF with history-wise dissipation in order to invoke Theorem \ref{prop-existing-LKF-expISS}. It crucially relies on the following technical lemma, which originally appeared in \cite[Lemma 5.3]{karafyllis2010necessary}.

\begin{lemma}[From point-wise to LKF-wise]\label{lemma_Iasson}
Given any continuously differentiable function $W:\mR^n\to\mRp$, the functional $V_0:\mathcal X^n\to\mRp$ defined as
\begin{align*}
V_0(\phi):=\max_{\tau\in[-\Delta,0]} e^{2\tau}W(\phi(\tau)),\quad \forall \phi\in\mathcal X^n,
\end{align*}
is Lipschitz on bounded sets and satisfies
\begin{align*}
V_0(\phi)>W(\phi(0)) \quad &\Rightarrow\quad D^+V_0(\phi,f(\phi,v))\leq -2 V_0(\phi) \\
V_0(\phi)= W(\phi(0))\quad &\Rightarrow\quad D^+V_0(\phi,f(\phi,v))\leq \max\left\{ -2V_0(\phi) ,\nabla W(\phi(0)) f(\phi,v)\right\}.
\end{align*}
\end{lemma}

We will see in the proof of Theorem \ref{theo-expISS2} that this result turns out to be particularly useful not only to construct a coercive LKF based on a non-coercive one, but also to impose a LKF-wise dissipation rate.

\subsection{Left growth condition}

The growth condition of Theorem \ref{theo-expISS2} imposes that solutions do not grow too fast. The following result shows that a similar statement can be obtained by requiring that they do not decay too fast. Its proof is provided in Section \ref{sec-proof-theo-expISS-left}.

\begin{theorem}[Exp-ISS under a left growth condition]\label{theo-expISS-left} Assume that there exist a functional $V:\mathcal X^n\to \mRp$ which is Lipschitz on bounded sets, $\underline a,\overline a,a>0$, $c\ge 0$ and $\gamma \in \mathcal N$ such that, for all $\phi \in \mathcal X^n$ and $v\in \mR^n$,
\begin{align} \label{GrindEQ__2_} 
\underline a \left|\phi (0)\right|^{2} &\le V(\phi )\le \overline a \left\| \phi \right\| ^{2}  \\ 
D^{+} V(\phi ,f(\phi ,v))&\le -a \left|\phi (0)\right|^{2} +c\left\| \phi \right\| ^{2} +\gamma \left(\left|v\right|\right).  \label{GrindEQ__3_} 
\end{align} 
Assume further that there exists a symmetric positive definite matrix $P\in \mR^{n\times n} $ and a constant $\sigma>0$ such that, for all $\phi \in \mathcal X^n$ and all $v\in \mR ^{m} $,
\begin{equation} \label{GrindEQ__4_} 
\phi(0)^\top P\, f(\phi ,v)\ge -\sigma\left(\left\| \phi \right\| ^{2} +\gamma \left(\left|v\right|\right)\right).
\end{equation} 
Then there exists $\bar c >0$ such that, if $c\in[0,\bar c)$, the system \eqref{eq-01} is exp-ISS. If, in addition, there exists $\gamma_0\geq 0$ such that $\gamma(s)=\gamma_0s^2$ for all $s\geq 0$, then \eqref{eq-01} is exp-ISS with linear gain.
\end{theorem}

Unlike Theorem \ref{theo-expISS2}, this result does require a lower bound on the LKF in terms of $|\phi(0)|$. Note however that $V$ is not requested to be coercive.

\begin{remark}[Estimate of $\bar c$]\label{rmk2}
Although Theorem \ref{theo-expISS-left} is an existence result, its proof provides an explicit estimate of the constant $\bar c$: see \eqref{eq-34} where $p_m,p_M  >0$ denote the minimal and maximal eigenvalues of $P$ respectively. In particular, $\bar c$ is independent of the function $\gamma$. However, here also the reader should be warned that the above estimate is conservative and in practice larger values of $c$ could be allowed. 
\end{remark}


The following result is an immediate consequence of Theorem \ref{theo-expISS-left} for input-free systems.

\begin{corollary}[GES under a left growth condition] Assume that there exist a functional $V:\mathcal X^n\to\mRp $ which is Lipschitz on bounded sets, $\underline a,\overline a, a>0$ and $c\geq 0$ such that, for all $\phi\in\mathcal X^n$,  \eqref{GrindEQ__2_} holds and
\begin{equation} \label{GrindEQ__6_} 
D^{+} V(\phi ,f_0(\phi ))\le -a  \left|\phi (0)\right|^{2} +c\left\| \phi \right\| ^{2}.
\end{equation} 
Assume further that there exists a symmetric positive definite matrix $P\in \mR^{n\times n} $ and a constant $\sigma>0$ such that, for all $\phi\in\mathcal X^n$,
\begin{equation} \label{GrindEQ__7_} 
\phi(0)^\top Pf_0(\phi )\ge -\sigma\left\| \phi \right\| ^{2}.
\end{equation} 
Then there exists a constant $\bar c >0$ such that, if $c\in[0,\bar c)$, the origin of \eqref{eq-0} is GES.
\end{corollary}

The proof of Theorem \ref{theo-expISS-left} relies on the following two lemmas, which may be of independent interest. The first one, proved in Section \ref{sec-proof-lem-halanay-RFC}, is used to show that, under the assumptions of Theorem \ref{theo-expISS-left}, the solutions of \eqref{eq-0} exist on $\mRp$. It actually ensures a stronger property called \emph{robust forward completeness (RFC)} \cite[Definition 2.1]{KAJIbook11} or \emph{bounded reachability property} \cite[Definition 4]{mironchenko2017characterizations}, which play a key role in the analysis of time delay systems. This property means that, for all $T,r>0$, all $x_0\in\mathcal X^n$, and all $u\in\mathcal U^m$, the corresponding solution of \eqref{eq-01} satisfies the following implication:
\begin{align*}
\|x_0\|\leq r,\ \|u\|\leq r\quad \Rightarrow\quad \sup_{t\in[0,T]} |x(t)|<+\infty.
\end{align*}
In other words, starting from any bounded set of initial states and considering any input lying in any given bounded set, solutions can only reach a bounded set over any finite time interval.

\begin{lemma}[Sufficient condition for RFC]\label{lem-halanay-RFC}
Assume that there exist a functional $V:\mathcal X^n\to\mRp$ which is Lipschitz on bounded sets, $a,c,\bar c\geq 0$, $\alpha,\overline\alpha\in\mathcal K_\infty$, and $\gamma \in\mathcal N$ such that, for all $\phi\in\mathcal X^n$ and all $v\in\mR^m$,
\begin{align}
\alpha(|\phi(0)|) &\leq V(\phi)\leq \overline \alpha(\|\phi\|)+\bar c\label{eq-40}\\
D^+V(\phi,f(\phi,v))&\leq a\alpha(\|\phi\|)+\gamma(|v|)+c.\label{eq-41}
\end{align}
Then the system \eqref{eq-01} is robustly forward complete.
\end{lemma}

From the regularity assumptions made on both $V$ and $f$, the upper bound \eqref{eq-41} is necessarily satisfied for some $\alpha\in\mathcal K_\infty$, $\gamma\in\mathcal N$, and $a,c\geq 0$. The key requirement in the above lemma therefore lies in the fact that $\alpha$ can also be picked as a lower bound on $V$, as imposed by \eqref{eq-40}.

The second lemma allows to convert the proof of exp-ISS into specific inequalities on the solutions norm.

\begin{lemma}[Exp-ISS characterization]\label{lem-expISS} The system \eqref{eq-01} is exp-ISS if and only if there exist constants $\ell,T>0$, $\lambda \in (0,1)$ and a function $\mu \in \mathcal N$ such that, for all $x_{0} \in \mathcal W^n$ and all $u\in\mathcal U^m$, its corresponding solution satisfies
\begin{align} 
\left\| x_{t} \right\| &\le \ell\left\| x_{0} \right\| +\mu \left(\left\| u_{[0,t]}\right\| \right), \quad \forall t\in [0,T],  \label{GrindEQ__9_}\\
\left\| x_{T} \right\| &\le \lambda \left\| x_{0} \right\| +\mu \left(\left\| u_{[0,T]}\right\| \right).  \label{GrindEQ__10_} 
\end{align} 
Moreover, \eqref{eq-01} is exp-ISS with linear gain if and only if the above holds with a linear function $\mu\in\mathcal N$.
\end{lemma}

Note that this result allows to consider only initial states $x_0$ that are absolutely continuous, which significantly increases the regularity of the resulting solutions. In particular it can be shown that, if $V$ is Lipschitz on bounded sets then, given any $x_{0} \in \mathcal W^n$ and any $u\in \mathcal U^m$, the map $t\mapsto V(x_t)$ is locally absolutely continuous (and consequently differentiable almost everywhere) on its domain of existence: see \cite[Theorem 5, Remark 6]{Pepe2007} or \cite[Lemma 2.5]{Karafyllis:2008hc}.

\section{Examples}

We illustrate our results through the following three academic examples.

\begin{example}[Theorems \ref{theo-expISS2} or \ref{theo-expISS-left}]\label{exa-1} 
Consider the system
\begin{subequations}\label{eq-exa1}
\begin{align}\label{eq-50}
\dot x_1(t) &= -\frac{1}{2}x_1(t)+x_2(t-\Delta)+x_2(t)\left(x_1(t)^2+x_2(t-\Delta)^2\right)\\
\dot x_2(t) &= -2x_2(t)-x_1(t)\left(x_1(t)^2+x_2(t-\Delta)^2\right)+u(t).
\end{align}
\end{subequations}
and let $f(x_t,u(t))$ denote its right-hand side. Consider the LKF defined for all $\phi=(\phi_1,\phi_2)^\top\in\mathcal X^2$ as
\begin{align}\label{eq-26}
V(\phi):=\phi_1(0)^2+\phi_2(0)^2+2\int_{-\Delta}^0\phi_2(\tau)^2d\tau.
\end{align}
Then it holds that
\begin{align}\label{eq-51}
|\phi(0)|^2\leq V(\phi)\leq (1+2\Delta)\|\phi\|^2,\quad \forall \phi\in\mathcal X^2,
\end{align}
which establishes \eqref{eq-2} and \eqref{GrindEQ__2_} with $\underline a=1$ and $\overline a=1+2\Delta$. For all $\phi\in\mathcal X^2$ and all $v\in\mR$, its Driver's derivative along the system's solutions reads 
\begin{align}
D^+V(\phi,f(\phi,v)) =&\, -\phi_1(0)^2+2\phi_1(0)\phi_2(-\Delta)-4\phi_2(0)^2+2\phi_2(0)v\nonumber\\
&\,+2\left(\phi_2(0)^2-\phi_2(-\Delta)^2\right)\nonumber\\
\leq&\, -\phi_1(0)^2+\frac{1}{2}\phi_1(0)^2+2\phi_2(-\Delta)^2-4\phi_2(0)^2+\phi_2(0)^2+v^2\nonumber\\
&\,+2\left(\phi_2(0)^2-\phi_2(-\Delta)^2\right)\nonumber\\
&\leq\,-\frac{1}{2}\phi_1(0)^2-\phi_2(0)^2+v^2\nonumber\\
&\leq\,-\frac{1}{2}|\phi(0)|^2+v^2,\label{eq-26bis}
\end{align}
thus establishing \eqref{eq-3} and \eqref{GrindEQ__3_}  with $a=1/2$ and $\gamma(s)=s^2$ for all $s\geq 0$.
Finally, 
\begin{align}
\phi(0)^\top f(\phi,v) &= -\frac{1}{2}\phi_1(0)^2+\phi_1(0)\phi_2(-\Delta)-2\phi_2(0)^2+\phi_2(0)v\nonumber\\
&\leq -\frac{1}{4}\phi_1(0)^2+\phi_2(-\Delta)^2-\frac{3}{2}\phi_2(0)^2+\frac{1}{2}v^2\nonumber\\
&\leq \,\|\phi\|^2+\gamma(|v|), \label{eq-28}
\end{align}
thus fulfilling \eqref{eq-22} with $\sigma=1$ and $P=I$. Hence, Theorem \ref{theo-expISS2} applies and we conclude that, for any value of the delay $\Delta\geq 0$, \eqref{eq-exa1} is exp-ISS with linear gain. Similarly,
\begin{align}
\phi(0)^\top f(\phi,v) &= -\frac{1}{2}\phi_1(0)^2+\phi_1(0)\phi_2(-\Delta)-2\phi_2(0)^2+\phi_2(0)v\nonumber\\
&\geq -\phi_1(0)^2-\frac{1}{2}\phi_2(-\Delta)^2-\frac{5}{2}\phi_2(0)^2-\frac{1}{2}v^2\nonumber\\
&\geq \,-3(\|\phi\|^2+\gamma(|v|)),
\end{align}
thus making the left growth condition \eqref{GrindEQ__4_} fulfilled with $\sigma=3$. So Theorem \ref{theo-expISS-left}. can also be invoked to conclude exp-ISS.

Note that the vector field in \eqref{eq-exa1} is not globally Lipschitz, thus making \cite[Corollary 1]{chaillet2019relaxed} inapplicable to conclude exp-ISS. Also, it is worth stressing that the considered LKF $V$ does not admit an LKF-wise dissipation rate, so Theorem \ref{theo-LKF-expISS-nogrowth} cannot be invoked either. It is however fair to stress that such a LKF-wise dissipation could have been obtained by considering the following LKF:
\begin{align*}
W(\phi)=\phi_1(0)^2+\phi_2(0)^2+\int_{-\Delta}^0 e^{-\rho\tau}\phi_2(\tau)^2d\tau,
\end{align*}
for some conveniently chosen $\rho>0$. The addition of an exponential term under the integral is standard in the analysis of time-delay systems. While this trick does work for some relevant classes of systems \cite[Lemma 1]{ORCHDESI21}, it is not universal and the possibility offered by Theorems \ref{theo-expISS2} and \ref{theo-expISS-left} to consider a point-wise dissipation significantly simplifies the analysis.

\end{example}

The second example illustrates how the $c\|\phi\|^2$ term appearing in \eqref{eq-3} and \eqref{GrindEQ__3_} can be used to guarantee some robustness to modeling errors.

\begin{example}[Theorem \ref{theo-expISS2} or \ref{theo-expISS-left} for robustness]
Consider the system
\begin{subequations}\label{eq-exa2bis}
\begin{align}
\dot x_1(t) &= -\frac{1}{2}x_1(t)+x_1(t-\Delta)+x_2(t)\left(x_1(t)^2+x_2(t-\Delta)^2\right)+\varepsilon d_1(x_t) \label{eq-exa21bis}\\
\dot x_1(t) &=-2x_2(t)-x_1(t)\left(x_1(t)^2+x_2(t-\Delta)^2\right)+u(t)+\varepsilon d_2(x_t).\label{eq-exa22bis}
\end{align}
\end{subequations}
The only difference with respect to \eqref{eq-exa1} lies in the extra terms $d_1(x_t)$ and $d_2(x_t)$ that represent modeling uncertainties. $\varepsilon\in\mR$ represents the intensity of these uncertainties. We assume that $d_1,d_2:\mathcal X^2\to\mR$ are Lipschitz on bounded sets and satisfy, for all $\phi\in\mathcal X^2$,
\begin{align}\label{eq-63}
|d_i(\phi)|\leq \|\phi\|, \quad \forall i\in\{1,2\}.
\end{align}
Considering the LKF $V$ defined in \eqref{eq-26}, it can be seen from \eqref{eq-26bis} and \eqref{eq-63} that
\begin{align*}
D^+V(\phi,f(\phi,v)) &\leq -\frac{1}{2}|\phi(0)|^2+v^2+2\varepsilon\left(\phi_1(0)d_1(\phi)+\phi_2(0)d_2(\phi)\right)\\
&\leq -\frac{1}{2}|\phi(0)|^2+v^2+4|\varepsilon|\|\phi\|^2,
\end{align*}
where $f$ is such that the right-hand side of \eqref{eq-exa2bis} reads $f(x_t,u(t))$. 
We conclude, with either Theorem \ref{theo-expISS2} or Theorem \ref{theo-expISS-left}, that there exists $\bar \varepsilon>0$ so that exp-ISS is preserved under such model uncertainties provided that $\varepsilon\in(-\bar\varepsilon,\bar\varepsilon)$. Recalling that, for this example, $\underline a=\sigma=p_m=p_M=1$ and $a=1/2$, the estimate of $\bar\varepsilon$ given by Theorem \ref{theo-expISS2} (see Remark \ref{rmk1}) reads $\bar\varepsilon_1 = e^{-4\Delta}/16$, thus suggesting a decrease of robustness when the delay $\Delta$ increases. Similarly, recalling that $\sigma=3$ in the application of Theorem \ref{theo-expISS-left}, we get from Remark \ref{rmk2} the following alternative estimate of this tolerable uncertainty:
\begin{align*}
\bar\varepsilon_2 = \frac{1}{8\left\lceil 2304(1+2\Delta)^4\right\rceil\left(\Delta+\frac{1}{48(1+2\Delta)^2}\right)},
\end{align*}
which can be checked to be greater than $\bar\varepsilon_1$ whenever $\Delta\geq 4.5$. Combining Theorems \ref{theo-expISS2} and \ref{theo-expISS-left} and letting $\bar\varepsilon:=\max\{\bar\varepsilon_1,\bar\varepsilon_2\}$, any $\varepsilon\in (-\bar \varepsilon,\bar\varepsilon)$ preserves exp-ISS of \eqref{eq-exa2bis} if the modeling uncertainty satisfies \eqref{eq-63}.
\end{example}

The following other variant of Example \ref{exa-1} shows that, on some occasions, condition \eqref{eq-22} can be fulfilled while condition \eqref{GrindEQ__4_} is violated.

\begin{example}[Theorem \ref{theo-expISS2}, but not Theorem \ref{theo-expISS-left}]\label{exa-3} 
Consider the system
\begin{subequations}\label{eq-exa2}
\begin{align}
\dot x_1(t) &= -\frac{1}{2}x_1(t)+x_1(t-\Delta)+x_2(t)\left(x_1(t)^2+x_2(t-\Delta)^2\right) \label{eq-exa21}\\
\dot x_1(t) &=-2x_2(t)-x_2(t)^3-x_1(t)\left(x_1(t)^2+x_2(t-\Delta)^2\right)+u(t).\label{eq-exa22}
\end{align}
\end{subequations}
The only difference with respect to \eqref{eq-exa1} is the extra term $-x_2(t)^3$ appearing in \eqref{eq-exa22}. This term brings only non-positive terms in Driver's derivative of the function $V$ defined in \eqref{eq-26}, meaning that \eqref{eq-26bis} remains valid when $f$ denotes the vector field in \eqref{eq-exa2}. For the same reason, \eqref{eq-28} remains valid as well. Theorem \ref{theo-expISS2} thus applies and the system \eqref{eq-exa2} is exp-ISS with linear gain. Nevertheless,
\begin{align*}
\phi(0)^\top f(\phi,v) &= -\frac{1}{2}\phi_1(0)^2+\phi_1(0)\phi_1(-\Delta)-2\phi_2(0)^2-\phi_2(0)^4+\phi_2(0)v.
\end{align*}
The term $-\phi_2(0)^4$ cannot be lower-bounded by a term proportional to $-\|\phi\|^2$. Reasoning similarly, we can see that \eqref{GrindEQ__4_}  cannot be fulfilled for this system, no matter the choice of the positive definite matrix $P\in\mR^{2\times 2}$, thus making Theorem \ref{theo-expISS-left} inapplicable.
\end{example}

 

%

\section{Proofs}

\subsection{Proof of Lemma \ref{lem-halanay-RFC}}\label{sec-proof-lem-halanay-RFC} 

Invoking the density of $\mathcal W^n$ in $\mathcal X^n$ and the continuity of solutions with respect to the initial state, it can be shown that robust forward completeness holds even if we restrict the analysis to initial states in $\mathcal W^n$. In other words, it is enough to show that, given any $r,T>0$, there exist $R>0$ such that, given any $x_0\in \mathcal W^n$ and any $u\in \mathcal U^m$, the following implication holds:
\begin{align}\label{eq-45}
\|x_0\|\leq r\ ,\|u\|\leq r \quad \Rightarrow\quad \sup_{t\in[0,T]} \left| x(t)\right| \leq R.
\end{align}
Consider any $x_0\in \mathcal W^n$ and any $u\in \mathcal U^m$ and let $x(\cdot)$ denote the maximal solution of \eqref{eq-01}, defined on an interval $[-\Delta,t_{\max})$ for some $t_{\max}\in[0,+\infty]$. From \cite[Theorem 2.1, p. 41]{Hale:1977tr}, we know that $t_{\max}\in(0,+\infty]$. The function $V$ being Lipschitz on bounded sets, it can be seen that the map $t\mapsto V(x_t)$ is then locally absolutely continuous on $[0,t_{\max})$: see \cite[Theorem 5, Remark 6]{Pepe2007} or \cite[Lemma 2.5]{Karafyllis:2008hc}. Hence, its derivative exists almost everywhere on $[0,t_{\max})$. Using \cite[Theorem 2]{pepe2007liapunov}, we get from \eqref{eq-41} that
\begin{align*}
\frac{d}{dt}V(x_t)\leq a\alpha(\|x_t\|)+\gamma(|u(t)|)+c,\quad \forall t\in[0,t_{\max})\ a.e.
\end{align*}
Integrating, it follows that, for all $t\in[0,t_{\max})$,
\begin{align*}
V(x_t)\leq V(x_0)+a\int_0^t \alpha(\|x_\tau\|) d\tau +\int_0^t \left(\gamma(|u(\tau)|)+c\right)d\tau.
\end{align*}
We then get from \eqref{eq-40} that, for all $t\in [0,t_{\max})$,
\begin{align}
\alpha(|x(t)|) \leq \overline \alpha(\|x_0\|)+\overline c+a\int_0^t \alpha(\|x_\tau\|) d\tau+\int_0^t \left(\gamma(|u(\tau)|)+c\right)d\tau. \label{eq-42}
\end{align}
It follows in particular that\footnote{Note that the interval $[\Delta,t_{\max})$ may be empty.}, for all $t\in [\Delta,t_{\max})$,
\begin{align}\label{eq-44}
\alpha(\|x_t\|) \leq\overline \alpha(\|x_0\|)+\overline c+a\int_0^t \alpha(\|x_\tau\|) d\tau+\int_0^t \left(\gamma(|u(\tau)|)+c\right)d\tau.
\end{align}
Moreover, for all $t\in [0,\min\{\Delta,t_{\max}\})$, it holds that
\begin{align*}
\alpha(\|x_t\|)  
&=\max\left\{\max_{\tau\in[t-\Delta,0]}\alpha(|x(\tau)|)\,,\, \max_{\tau\in[0,t]}\alpha(|x(\tau)|) \right\}\\
&\leq \max\left\{\alpha(\|x_0\|)\,,\, \max_{\tau\in[0,t]}\alpha(|x(\tau)|) \right\}.
\end{align*}
Using again \eqref{eq-42} and observing that \eqref{GrindEQ__2_} imposes that $\overline \alpha(s)+\overline c\geq \alpha(s)$ for all $s\geq 0$, it follows that, for all $t\in [0,\min\{\Delta,t_{\max}\})$,
\begin{align}
\alpha(\|x_t\|)  &\leq \max\left\{\alpha(\|x_0\|) \,,\, \overline \alpha(\|x_0\|)+\overline c+a\int_0^t \alpha(\|x_\tau\|) d\tau+\int_0^t \left(\gamma(|u(\tau)|)+c\right)d\tau\right\}\nonumber\\
&\leq\overline \alpha(\|x_0\|)+\overline c+a\int_0^t \alpha(\|x_\tau\|) d\tau+\int_0^t \left(\gamma(|u(\tau)|)+c\right)d\tau.\label{eq-43}
\end{align}
In view of \eqref{eq-44}, we conclude that \eqref{eq-43} holds for all $t\in[0,t_{\max})$. Invoking Gronwall-Bellman's inequality, it follows that, for all $t\in [0,t_{\max})$,
\begin{align*}
\alpha(\|x_t\|)  \leq \left(\overline \alpha(\|x_0\|)+\overline c\right)e^{at}+\int_0^t \left(\gamma(|u(\tau)|)+c\right)e^{a(t-\tau)}d\tau,
\end{align*}
meaning that
\begin{align}\label{eq-46}
\|x_t\|  \leq \alpha^{-1}\left(\left(\overline \alpha(\|x_0\|)+\overline c\right)e^{at}+\int_0^t \left(\gamma(|u(\tau)|)+c\right)e^{a(t-\tau)}d\tau\right).
\end{align}
If $t_{\max}$ were finite, then we would have $\limsup_{t\to t_{\max}^-}\|x_t\|=+\infty$ \cite[Theorem 3.2, p. 43]{Hale:1977tr}, which is incompatible with \eqref{eq-46}. Hence $t_{\max}=+\infty$. In particular, given any $r,T>0$, if $\|x_0\|\leq r$ and $\|u\|\leq r$, then it holds that
\begin{align*}
\sup_{t\in[0,T]} |x(t)|\leq \alpha^{-1}\left(\left(\overline \alpha(r)+\overline c\right)e^{aT}+(\gamma(r)+c)\frac{e^{aT}}{a}
\right),
\end{align*}
thus fulfilling \eqref{eq-45} and establishing robust forward completeness.

\subsection{Proof of Theorem \ref{theo-LKF-expISS-nogrowth}}\label{sec-proof-theo-LKF-expISS-nogrowth}

First observe that Lemma \ref{lem-halanay-RFC} readily ensures that all solutions of \eqref{eq-01} exist on $\mRp$. Using density of $C^{1} \left([-\Delta,0];{\mathbb R}^{n} \right)$ in $\mathcal X^n$ and continuity of solutions with respect to the initial state, it is enough to show that there exist $k,\eta>0$ and $\mu\in\mathcal N$ such that, for all $x_0\in C^{1} \left([-\Delta,0];{\mathbb R}^{n} \right)$ and all $u\in\mathcal U^m$, 
\begin{align}\label{eq-62}
|x(t)|\leq k\|x_0\|e^{-\eta t} + \mu(\|u_{[0,t]}\|),\quad \forall t\geq 0.
\end{align}
Considering any $x_{0} \in C^{1} \left([-\Delta,0];{\mathbb R}^{n} \right)$ and any locally bounded $u\in \mathcal U^m$, since $V$ is Lipschitz on bounded sets, the mapping $t\mapsto V\left(x_{t} \right)$ is locally absolutely continuous on $\mRp$ (see \cite[Theorem 5, Remark 6]{Pepe2007} or \cite[Lemma 2.5]{Karafyllis:2008hc}) and we get from \eqref{eq-58} that
\begin{equation*} 
\frac{d}{dt} V\left(x_{t} \right)\le -a V\left(x_{t} \right)+c\left\| x_{t} \right\| ^{\rho} +\gamma \left(\left|u(t)\right|\right),\quad \forall t\geq 0\ a.e.
\end{equation*} 
It follows that
\begin{equation*} 
\frac{d}{dt} \left(e^{at}V\left(x_{t} \right)\right)\le ce^{at}\left\| x_{t} \right\| ^{\rho} +e^{at}\gamma \left(\left|u(t)\right|\right),\quad \forall t\geq 0\ a.e.
\end{equation*} 
Since $t\mapsto V\left(x_{t} \right)$ is locally absolutely continuous, we get for all $t\geq 0$ that
\begin{align*} 
e^{at}V\left(x_{t} \right)-V\left(x_{0} \right)&=\int _{0}^{t}\frac{d}{d\tau} \left(e^{a\tau}V(x_\tau)\right)d\tau \\ 
&\le c\int _{0}^{t}e^{a\tau}\left\| x_{\tau} \right\| ^{\rho} d\tau +\int _{0}^{t}e^{a\tau}\gamma \left(\left|u(\tau)\right|\right)d\tau.
\end{align*} 
It follows from \eqref{eq-57} that, for all $t\geq 0$,
\begin{equation} 
\left|x(t)\right|^{\rho} \le \frac{e^{-at}}{\underline a}\left(\overline a \left\| x_{0} \right\| ^{\rho} +c\int _{0}^{t}e^{a\tau}\left\| x_{\tau} \right\| ^{\rho} d\tau +\int _{0}^{t}e^{a\tau}\gamma \left(\left|u(\tau)\right|\right)d\tau \right). \label{eq-59}
\end{equation} 
In particular, for all $t\geq \Delta$, we get that
\begin{align}\label{eq-59bis}
\|x_t\|^{\rho} \le \frac{e^{-a(t-\Delta)}}{\underline a}\left(\overline a \left\| x_{0} \right\| ^{\rho} +c\int _{0}^{t}e^{a\tau}\left\| x_{\tau} \right\| ^{\rho} d\tau +\int _{0}^{t}e^{a\tau}\gamma \left(\left|u(\tau)\right|\right)d\tau \right). 
\end{align}
Moreover, using again \eqref{eq-59}, it holds  for all $t\in [0,\Delta]$ that
\begin{align*}
\|x_t\|^\rho&=\max\left\{\max_{\tau\in[t-\Delta,0]} |x(\tau)|^\rho\,,\,\max_{\tau\in[0,t]} |x(\tau)|^\rho\right\}\\
&\leq \max\left\{\|x_0\|^\rho\,,\,\frac{1}{\underline a}\left(\overline a \left\| x_{0} \right\| ^{\rho} +c\int _{0}^{t}e^{a\tau}\left\| x_{\tau} \right\| ^{\rho} d\tau +\int _{0}^{t}e^{a\tau}\gamma \left(\left|u(\tau)\right|\right)d\tau \right)\right\}\\
&\leq\frac{1}{\underline a}\left(\overline a \left\| x_{0} \right\| ^{\rho} +c\int _{0}^{t}e^{a\tau}\left\| x_{\tau} \right\| ^{\rho} d\tau +\int _{0}^{t}e^{a\tau}\gamma \left(\left|u(\tau)\right|\right)d\tau \right),
\end{align*}
where the last bound comes from the fact that $\overline a\geq \underline a$ as ensured by \eqref{eq-57}. We conclude that \eqref{eq-59bis} actually holds for all $t\geq 0$, which implies that, for all $t\geq 0$,
\begin{align*} 
\left\| x_{t} \right\| ^{\rho} \le \frac{\overline a e^{-a(t-\Delta)}}{\underline a}\left\| x_{0} \right\|^{\rho} +\frac{e^{a\Delta}}{\underline a a}\sup_{\tau\in[0,t]} \left\{c\left\| x_{\tau} \right\| ^{\rho} +\gamma \left(\|u_{[0,\tau]}\|\right)\right\}.
\end{align*} 
Invoking \cite[Lemma 7.1]{KAKR19book}, we conclude that for every $\varepsilon >0$ there exists $\delta >0$ (independent of the solution) such that, for all $t\ge 0$,
\begin{align*}
\left\| x_{t} \right\| ^{\rho} \le&\, \frac{\overline ae^{a\Delta}}{\underline a}\|x_0\|^\rho e^{-\delta t} +\frac{e^{a\Delta}}{\underline a a} (1+\varepsilon )\sup_{\tau\in[0,t]}\left\{\left(c\left\| x_{\tau} \right\| ^{\rho} +\gamma(\|u_{[0,\tau]}\|)\right)e^{-\delta(t-\tau)}\right\},
\end{align*} 
which in turn implies that
\begin{align} \label{GrindEQ__51_} 
\sup_{\tau\in[0,t]} \left\| x_{\tau} \right\| ^{\rho} e^{\delta\tau}\le&\, \frac{\overline a e^{a\Delta}}{\underline a}\left\| x_{0} \right\| ^{\rho}  +\frac{ce^{a\Delta}}{\underline a a} (1+\varepsilon )\sup_{\tau\in[0,t]} \left\{\left\| x_{\tau} \right\| ^{\rho}e^{\delta\tau}\right\}\nonumber \\ 
&\,+\frac{e^{a\Delta}}{\underline a a} (1+\varepsilon )\gamma(\|u_{[0,t]}\|)e^{\delta t}.
\end{align} 
Pick $\varepsilon >0$ small enough that $
ce^{a\Delta} (1+\varepsilon )/\underline a a<1$,
the existence of which is ensured by \eqref{eq-56}. 
Then the constant $\xi:=1-ce^{a\Delta} (1+\varepsilon )/\underline a a$ is positive and we get from \eqref{GrindEQ__51_}  that, for all $t\geq 0$,
\begin{align*} 
\|x_t\|^\rho e^{\delta t}\leq \sup_{\tau\in[0,t]}\left\| x_{\tau} \right\| ^{\rho} e^{\delta\tau}\le \frac{\overline a e^{a\Delta}}{\underline a\xi}\left\| x_{0} \right\| ^{\rho}  +\frac{e^{a\Delta}}{\underline a a\xi} (1+\varepsilon )\gamma(\|u_{[0,t]}\|)e^{\delta t}.
\end{align*} 
The function $s\mapsto s^{1/\rho}$ being increasing, it holds that $(r+s)^{1/\rho}\leq (2r)^{1/\rho}+(2s)^{1/\rho}$ for all $r,s\geq 0$. We conclude that
\begin{align*} 
\|x_t\| &\leq \left(\frac{2\overline a e^{a\Delta}}{\underline a\xi}\right)^{1/\rho}\left\| x_{0} \right\|e^{-\delta t/\rho}  +\left(\frac{2e^{a\Delta}}{\underline a a\xi} (1+\varepsilon )\right)^{1/\rho}\gamma \left(\|u_{[0,t]}\|\right)^{1/\rho},
\end{align*} 
which establishes \eqref{eq-62} with $\eta:=\delta$ and
\begin{align*}
k:=\left(\frac{2\overline a e^{a\Delta}}{\underline a\xi}\right)^{1/\rho},\quad
\mu(s):=\left(\frac{2e^{a\Delta}}{\underline a a\xi} (1+\varepsilon )\right)^{1/\rho}\gamma \left(s\right)^{1/\rho},\quad \forall s\geq 0.
\end{align*}
We conclude that \eqref{eq-01} is exp-ISS. Finally, if $\gamma(s)=\gamma_0s^{\rho}$ for some $\gamma_0\geq 0$, then the obtained gain $\mu$ is linear and \eqref{eq-01} is exp-ISS with linear gain.

\subsection{Proof of Theorem \ref{theo-expISS2}}\label{sec-proof-theo-expISS2}

Consider the functional $V_0:\mathcal X^n\to\mathbb R_{\geq 0}$ defined as
\begin{align}\label{eq-11}
V_0(\phi):=\max_{\tau\in[-\Delta,0]} e^{2\tau}\phi(\tau)^\top P\phi(\tau),\quad \forall \phi\in\mathcal X^n.
\end{align}
Then it holds that
\begin{align}\label{eq-14}
e^{-2\Delta}p_m\|\phi\|^2\leq V_0(\phi)\leq p_M\|\phi\|^2,
\end{align}
where $p_m,p_M>0$ denote the smallest and largest eigenvalues of $P$ respectively.
In view of Lemma \ref{lemma_Iasson}, $V_0$ is also Lipschitz on bounded sets and satisfies
\begin{align}
V_0(\phi)>\phi(0)^\top P\phi(0) \quad &\Rightarrow\quad D^+V_0(\phi,f(\phi,v))\leq -2 V_0(\phi) \label{eq-5}\\
V_0(\phi)= \phi(0)^\top P\phi(0)\quad &\Rightarrow\quad D^+V_0(\phi,f(\phi,v))\leq 2\max\left\{ -V_0(\phi) ,\phi(0)^\top P f(\phi,v)\right\}.\label{eq-6}
\end{align}
Now let $W(\phi):=V(\phi)+\varepsilon V_0(\phi)$ for some $\varepsilon>0$ to be chosen later. Then $W$ is Lipschitz on bounded sets and it holds that
\begin{align}\label{eq-55}
\varepsilon e^{-2\Delta}p_m\|\phi\|^2\leq W(\phi)\leq (\overline a+\varepsilon p_M)\|\phi\|^2.
\end{align}
Thus, $W$ is a coercive LKF. In view of \eqref{eq-5}-\eqref{eq-6}, we need to consider two cases.

\noindent\emph{Case 1: $V_0(\phi)>\phi(0)^\top P\phi(0)$.} Then we get from \eqref{eq-3} and \eqref{eq-5} that 
\begin{align}
D^+W(\phi,f(\phi,v))&\leq -2\varepsilon V_0(\phi)-a|\phi(0)|^2+c\|\phi\|^2+\gamma(|v|)\nonumber
\\&\leq-\left(2\varepsilon p_m e^{-2\Delta}-c\right)\|\phi\|^2+\gamma(|v|), \label{eq-9}
\end{align}
where the last bound comes from \eqref{eq-14}.

\noindent\emph{Case 2: $V_0(\phi)=\phi(0)^\top P\phi(0)$.} Then we have from \eqref{eq-14} that 
\begin{align}
|\phi(0)|^2\geq \frac{V_0(\phi)}{p_M}\geq \frac{p_m}{p_M}e^{-2\Delta}\|\phi\|^2.
\end{align}
It follows from \eqref{eq-3} and \eqref{eq-6} that
\begin{align}
D^+W(\phi,f(\phi,v)) \leq&\, -a|\phi(0)|^2+c\|\phi\|^2+\gamma(|v|)\nonumber\\
&\,+2\varepsilon\max\left\{-V_0(\phi)\,,\,\phi(0)^\top Pf(\phi,v)\right\}\nonumber\\
\leq&\, -\left(a\frac{p_m}{p_M}e^{-2\Delta}-c\right)\|\phi\|^2+\gamma(|v|)\nonumber\\
&\,+2\varepsilon\max\left\{-V_0(\phi)\,,\,\phi(0)^\top Pf(\phi,v)\right\}. \label{eq-24}
\end{align}
Here again, we need to consider two subcases.

\noindent\emph{Subcase 2.1: $-V_0(\phi)\geq \phi(0)^\top Pf(\phi,v)$.} Then we get from \eqref{eq-24} that 
\begin{align}
D^+W(\phi,f(\phi,v)) &\leq\,-\left(a\frac{p_m}{p_M}e^{-2\Delta}-c\right)\|\phi\|^2+\gamma(|v|)-2\varepsilon V_0(\phi)\nonumber\\
&\leq\,-\left(a\frac{p_m}{p_M}e^{-2\Delta}-c\right)\|\phi\|^2+\gamma(|v|).\label{eq-12}
\end{align}

\noindent\emph{Subcase 2.2: $-V_0(\phi) <\phi(0)^\top Pf(\phi,v)$.} Then we make use of condition \eqref{eq-22} to get from \eqref{eq-24} that
\begin{align*}
D^+W(\phi,f(\phi,v)) \leq &\,-\left(a\frac{p_m}{p_M}e^{-2\Delta}-c\right)\|\phi\|^2+\gamma(|v|) +2\varepsilon\phi(0)^\top Pf(\phi,v)\nonumber\\
\leq &\,-\left(a\frac{p_m}{p_M}e^{-2\Delta}-c\right)\|\phi\|^2+\gamma(|v|)+2\varepsilon\sigma(\|\phi\|^2+\gamma(|v|))\\
\leq &\,-\left(a\frac{p_m}{p_M}e^{-2\Delta}-c-2\varepsilon\sigma\right)\|\phi\|^2+(1+2\varepsilon\sigma)\gamma(|v|).
\end{align*}
By picking
\begin{align}\label{eq-30}
\varepsilon:= a\frac{p_m}{4\sigma p_M}e^{-2\Delta},
\end{align}
it follows in particular that
\begin{align}
D^+W(\phi,f(\phi,v)) &\leq -\left(a\frac{p_m}{2p_M}e^{-2\Delta}-c\right)\|\phi\|^2+(1+2\varepsilon\sigma)\gamma(|v|).
\label{eq-13}
\end{align}
Combining \eqref{eq-9}, \eqref{eq-12} and \eqref{eq-13}, we obtain the following inequality, valid for all $\phi\in\mathcal X^n$ and all $v\in\mathbb R^m$:
\begin{align*}
D^+W(\phi,f(\phi,v))\leq&\, -(\bar c-c)\|\phi\|^2+\bar \gamma(|v|)
\end{align*}
where
\begin{align}
\bar c&:=\min\left\{ 2\varepsilon\,,\,\frac{a}{2p_M}\right\}p_me^{-2\Delta}\label{eq-29} \\
\bar\gamma(s)&:=(1+2\varepsilon\sigma)\gamma(s),\quad \forall s\geq 0.\label{eq-25}
\end{align}
Using \eqref{eq-55}, we obtain that
\begin{align*}
D^+W(\phi,f(\phi,v))\leq&\, -\frac{\bar c-c}{\overline a+\varepsilon p_M}W(\phi)+\bar \gamma(|v|).
\end{align*}
We conclude that $W$ has a history-wise dissipation rate provided that $c\in[0,\bar c)$. Exp-ISS of \eqref{eq-01} then follows from classical LKF sufficient conditions for exp-ISS, such as the one recalled as Proposition \ref{prop-existing-LKF-expISS}. In the case when $\gamma$ is a square function, we see from \eqref{eq-25} that $\bar\gamma$ results in a square function. Exp-ISS with linear gain then also follows from Proposition \ref{prop-existing-LKF-expISS}.

\subsection{Proof of Lemma \ref{lem-expISS}}

If \eqref{eq-01} is exp-ISS then there exist constants $k\geq 1$, $\eta>0$ and a function $\mu \in \mathcal N$ such that, for all $x_{0} \in \mathcal W^n$ and all $u\in \mathcal U^m$, its solution satisfies
\begin{align}\label{eq-64}
|x(t)|\leq k\|x_0\|e^{-\eta t}+\mu(\|u_{[0,t]}\|),\quad \forall t\geq 0,
\end{align} 
which implies in particular that
\begin{align}\label{eq-65}
\|x_t\|\leq k\|x_0\|e^{-\eta t}+\mu(\|u_{[0,t]}\|),\quad \forall t\geq \Delta.
\end{align} 
In addition, we have from \eqref{eq-64} that, for all $t\in[0,\Delta]$,
\begin{align*}
\|x_t\| &= \max\left\{ \max_{\tau\in[t-\Delta,0]} |x(\tau)|\,,\,  \max_{\tau\in[0,t]} |x(\tau)|\right\}\\
&\leq \max\left\{ \|x_0\|\,,\,  k\|x_0\|e^{-\eta t}+\mu(\|u_{[0,t]}\|)\right\}.
\end{align*}
Since $k\geq 1$, it follows that, for all $t\in[0,\Delta]$,
\begin{align}\label{eq-66}
\|x_t\| \leq k\|x_0\|e^{-\eta (t-\Delta)}+\mu(\|u_{[0,t]}\|),\quad \forall t\in[0,\Delta].
\end{align}
In view of \eqref{eq-65}, we conclude that \eqref{eq-66} actually holds for all $t\geq 0$. In particular, given any $\lambda \in (0,1)$, inequalities \eqref{GrindEQ__9_} and \eqref{GrindEQ__10_} hold with $\ell:=ke^{\eta\Delta}$ and $T:=\Delta+\frac{1}{\eta } \ln \left(\frac{k}{\lambda } \right)$. 

We now proceed to showing the converse. To that aim, suppose that there exist $\ell,T>0$, $\lambda \in (0,1)$ and $\mu \in \mathcal N$ such that, for all $x_{0} \in \mathcal W^n$ and all $u\in \mathcal U^m$, the solution satisfies \eqref{GrindEQ__9_} and \eqref{GrindEQ__10_}. Invoking the density of $\mathcal W^n$ in $\mathcal X^n$ and the continuity of solutions with respect to the initial state, it can be shown that \eqref{GrindEQ__9_} and \eqref{GrindEQ__10_} actually hold for all $x_{0} \in \mathcal X^n$ and all  $u\in \mathcal U^m$. Considering any $x_{0} \in \mathcal X^n$ and any $u\in \mathcal U^m$, let $x(\cdot)$ denote the corresponding solution of \eqref{eq-01}. Clearly, \eqref{GrindEQ__9_}-\eqref{GrindEQ__10_} and the semigroup property imply that $x(\cdot)$ is defined on $\mRp$. Using \eqref{GrindEQ__10_} inductively in conjunction with the semigroup property, we conclude that, for every $i\in\mathbb N_{\ge 1}$, 
\begin{align*} 
\left\| x_{iT} \right\| \le \lambda ^{i} \left\| x_{0} \right\| +\mu\left(\left\| u_{[0,iT]}\right\| \right)\sum _{j=0}^{i-1}\lambda ^{j}\leq \lambda ^{i} \left\| x_{0} \right\| +\frac{1}{1-\lambda}\mu\left(\left\| u_{[0,iT]}\right\| \right).
\end{align*} 
Defining $\eta :=\frac{1}{T} \ln \left(\frac{1}{\lambda } \right)>0$,
we get that
\begin{equation*} 
\left\| x_{iT} \right\| \le e^{-i\eta T}\left\| x_{0} \right\| +\frac{1}{1-\lambda } \mu \left(\left\| u_{[0,iT]}\right\| \right).
\end{equation*} 
Using \eqref{GrindEQ__9_} and the semigroup property, it follows that, for every $i\in\mathbb N_{\ge 1}$ and for every $\tau\in [0,T]$,
\begin{equation*} 
\left\| x_{iT+\tau} \right\| \le \ell e^{-i\eta T}\left\| x_{0} \right\| +\left(\frac{\ell}{1-\lambda } +1\right)\mu \left(\left\| u_{[0,iT+\tau]}\right\| \right).
\end{equation*} 
Notice that, by \eqref{GrindEQ__9_}, the above estimate holds for $i=0$ as well.  Since for every $t\ge 0$ there exists $i\in\mathbb N$ and $\tau\in [0,T)$ such that $t=iT+\tau$, we obtain that
\begin{align*} 
\left\| x_{t} \right\|
&\le \ell e^{-\eta(t-\tau)}\left\| x_{0} \right\| +\left(\frac{\ell}{1-\lambda } +1\right)\mu \left(\left\| u_{[0,t]}\right\| \right)\\ 
&\le \ell e^{\eta T}\left\| x_{0} \right\|e^{-\eta t} +\left(\frac{\ell}{1-\lambda } +1\right)\mu \left(\left\| u_{[0,t]}\right\| \right),
\end{align*} 
and exp-ISS follows. If $\mu$ is linear then the obtained ISS gain is linear too.

\subsection{Proof of Theorem \ref{theo-expISS-left}}\label{sec-proof-theo-expISS-left}

Let $p_m,p_M>0$ denote respectively the smallest and largest eigenvalues of $P$. We assume that \eqref{GrindEQ__2_}, \eqref{GrindEQ__3_} and \eqref{GrindEQ__4_} hold for some $c\in \left(0,\bar c \right]$, where $\bar c$ is a positive constant to be defined later. The proof for the case $c=0$ is completely similar (and far easier). First observe that, by Lemma \ref{lem-halanay-RFC}, the system is robustly forward complete. In particular, its solutions exist on $\mRp$. The proof relies on Lemma \ref{lem-expISS}. Given any $x_{0} \in \mathcal W^n$ and any $u\in \mathcal U^m$, let $x(\cdot)$ denote the corresponding solution of  \eqref{eq-01}. Since $V$ is Lipschitz on bounded sets, the mapping $t\mapsto V(x_t)$ is then locally absolutely continuous (as ensured by \cite[Theorem 5, Remark 6]{Pepe2007} or \cite[Lemma 2.5]{Karafyllis:2008hc}), and we get from \eqref{GrindEQ__3_} that
\begin{equation*} 
\frac{d}{d\, t} V\left(x_{t} \right)\le -a  \left|x(t)\right|^{2} +c\left\| x_{t} \right\| ^{2}+\gamma \left(\left|u(t)\right|\right),\quad \forall t\geq 0\ a.e.
\end{equation*} 
Integrating, we get that, for all $t\geq 0$,
\begin{equation} \label{GrindEQ__19_} 
V\left(x_{t} \right)+a  \int _{0}^{t}\left|x(\tau)\right|^{2} d\tau \le V\left(x_{0} \right)+ct\max_{\tau\in[0,t]} \left\| x_{\tau} \right\| ^{2}+t\gamma \left(\left\| u_{[0,t]}\right\| \right).
\end{equation} 
Using \eqref{GrindEQ__2_}, it follows in particular that, for all $t\geq 0$,
\begin{align*} 
\left|x(t)\right|^{2}& \le \frac{\overline a }{\underline a } \left\| x_{0} \right\| ^{2} +\frac{ct}{\underline a } \max_{\tau\in[-\Delta,t]}\left|x(\tau)\right|^{2}+\frac{t}{\underline a } \gamma \left(\left\| u_{[0,t]}\right\| \right).
\end{align*} 
Consequently, for all $t\geq 0$,
\begin{align*} 
\max_{\tau\in[-\Delta,t]} |x(\tau)|^2 &=\max\left\{ \max_{\tau\in[-\Delta,0]} |x(\tau)|^2\,,\,\max_{\tau\in[0,t]} |x(\tau)|^2\right\}\\
&\leq\max\left\{ \|x_0\|^2\,,\,\frac{\overline a }{\underline a } \left\| x_{0} \right\| ^{2} +\frac{ct}{\underline a } \max_{\tau\in[-\Delta,t]}\left|x(\tau)\right|^{2}+\frac{t}{\underline a } \gamma \left(\left\| u_{[0,t]}\right\| \right)\right\}\\
&\leq \frac{\overline a }{\underline a } \left\| x_{0} \right\| ^{2} +\frac{ct}{\underline a } \max_{\tau\in[-\Delta,t]}\left|x(\tau)\right|^{2}+\frac{t}{\underline a } \gamma \left(\left\| u_{[0,t]}\right\| \right),
\end{align*} 
where we used the fact that $\overline a \geq \underline a$ as ensured by \eqref{GrindEQ__2_}. In particular, for all $t\in [0,\frac{\underline a }{2c}]$, it holds that
\begin{align}  
\max_{\tau\in[-\Delta,t]}|x(\tau)|^2
& \le \frac{\overline a}{\underline a-ct}\left\| x_{0} \right\| ^{2} +\frac{t}{\underline a-ct}\gamma \left(\left\| u_{[0,t]}\right\| \right)\nonumber\\
&\le \frac{2\overline a}{\underline a}\left\| x_{0} \right\| ^{2} +\frac{2t}{\underline a}\gamma \left(\left\| u_{[0,t]}\right\| \right),\label{GrindEQ__21_} 
\end{align} 
which in turn implies that, for all $T\in [0,\frac{\underline a }{2c}]$ and all $t\in[0,T]$,
\begin{align} \label{GrindEQ__23_} 
\left\| x_{t} \right\| &\le\sqrt{\frac{2\overline a}{\underline a}} \left\| x_{0} \right\|  +\sqrt{\frac{2T}{\underline a} }\sqrt{\gamma \left(\left\| u_{[0,t]}\right\| \right)} .
\end{align} 
On the other hand, we have from \eqref{GrindEQ__2_} and \eqref{GrindEQ__19_} that, for all $t\in[0,T]$, 
\begin{equation} \label{GrindEQ__25_} 
\int _{0}^{T}\left|x(\tau)\right|^{2} d\tau \le \frac{\overline a }{a  } \left\| x_{0} \right\| ^{2} +\frac{cT}{a  } \max_{\tau\in[0,t]}\left\| x_{\tau} \right\| ^{2}+\frac{T}{a  } \gamma \left(\left\| u_{[0,T]}\right\| \right).
\end{equation} 
Combining \eqref{GrindEQ__21_} and \eqref{GrindEQ__25_}, we get for all $T\in \left[0,\frac{\underline a }{2c} \right]$: 
\begin{align} 
\int _{0}^{T}\left|x(\tau)\right|^{2} d\tau &\le \frac{\overline a}{a}\left(1+\frac{2cT}{\underline a}\right)\|x_0\|^2+\frac{T}{a}\left(1+\frac{2cT}{\underline a}\right)\gamma(\|u_{[0,T]}\|) \nonumber\\
&\le \frac{2\overline a}{a}\|x_0\|^2+\frac{2T}{a}\gamma(\|u_{[0,T]}\|). \label{GrindEQ__26_} 
\end{align} 
Pick
\begin{align}\label{eq-34}
\varepsilon:=\frac{p_m\underline a^2}{16\overline a^2\sigma} \quad
q:=\left\lceil\frac{p_M}{\sigma \varepsilon^2}\right\rceil, \quad 
T := q(\Delta+\varepsilon),\quad 
\bar c:=\frac{\underline a}{2T},
\end{align}
and consider the function $\mu\in\mathcal N$ defined for all $s\geq 0$ as
\begin{align}
\mu(s)&:= \frac{4}{\sqrt{\underline a}}\max\left\{\sqrt{2T}\,,\, 
\sqrt{\frac{\overline a}{p_m}\left(\frac{p_MT}{aq\varepsilon}+\sigma\varepsilon\left(1+\frac{2T}{\underline a}\right) \right)}\right\}\sqrt{\gamma(s)}. \label{eq-38}
\end{align}
Then it follows readily from the definition of $\bar c$ in \eqref{eq-34} that, for all $c\in(0,\bar c]$,
\begin{align}\label{eq-32}
T\in\left(0,\frac{\underline a}{2c}\right], 
\end{align}
and it can be checked from the definitions of $q$ and $\varepsilon$ in \eqref{eq-34} that
\begin{align*}
\frac{q\varepsilon}{p_M}\left(\frac{p_m\underline a}{2\overline a}-\frac{4\varepsilon\sigma\overline a}{\underline a}\right)>\frac{2\overline a}{\underline a}.
\end{align*}
Consequently, there exists $\lambda\in(0,1)$ such that 
\begin{align}\label{eq-33}
\frac{q\varepsilon}{p_M}\left(\frac{p_m\underline a}{2\overline a}\lambda^2-\frac{4\varepsilon\sigma\overline a}{\underline a}\right)>\frac{2\overline a}{\underline a}.
\end{align}
Also, considering separately the two possible values that the $\max$ in \eqref{eq-38} may take, it can be checked that the following two properties hold for all $s\geq 0$:
\begin{align}
&\mu(s)-\sqrt{\frac{2T}{\underline a}}\sqrt{\gamma(s)}\geq \frac{\mu(s)}{2}\label{eq-39bis}\\
&\frac{q\varepsilon}{p_M}\left(\frac{\underline a p_m}{8\overline a}\mu(s)^2-2\sigma\varepsilon\left(1+\frac{2T}{\underline a}\right)\gamma(s)\right)\geq \frac{2T}{a}\gamma(s).\label{eq-39}
\end{align}
We show next that the estimate 
\begin{align}\label{eq-***}
\left\| x_{T} \right\| \le \lambda \left\| x_{0} \right\| +\mu\left(\left\| u_{[0,T]}\right\| \right)                                            
\end{align}
holds with the constant $T$ given in \eqref{eq-34}, the function $\mu$ provided in \eqref{eq-38}, and any $\lambda\in(0,1)$ satisfying \eqref{eq-33} with the constants provided in \eqref{eq-34}. The proof of \eqref{eq-***} is made by contradiction. Suppose that it does not hold for some $x_{0} \in \mathcal W^n$ and some $u\in \mathcal U^m$, namely: 
\begin{align}\label{eq-48}
\left\| x_{T} \right\| >\lambda \left\| x_{0} \right\| +\mu(\left\| u_{[0,T]}\right\|).
\end{align}
We know from \eqref{GrindEQ__23_} and the semi-group property that
\begin{align*}
\left\| x_{T} \right\| \le \sqrt{\frac{2\overline a}{\underline a} }\left\| x_{t} \right\|  +\sqrt{\frac{2T}{\underline a}}\sqrt{\gamma \left(\left\| u_{[0,T]}\right\| \right)} , \quad \forall t\in [0,T].
\end{align*}
This, combined with \eqref{eq-48}, gives for all $t\in[0,T]$:
 \begin{align*} 
\left\| x_{t} \right\| &>\sqrt{\frac{\underline a}{2\overline a}}\left(\lambda\left\| x_{0} \right\| +\mu(\|u_{[0,T]}\|)-\sqrt{\frac{2T}{\overline a}}\sqrt{\gamma(\|u_{[0,T]}\|)}\right).
\end{align*} 
It follows that, for each $i\in\{0,\ldots,q-1\}$, there exists $t_{i} \in \left[i(\Delta+\varepsilon ),i(\Delta+\varepsilon )+\Delta\right]$ such that 
\begin{align}
|x(t_i)| &>\sqrt{\frac{\underline a}{2\overline a}}\left(\lambda\left\| x_{0} \right\| +\mu(\|u_{[0,T]}\|)-\sqrt{\frac{2T}{\overline a}}\sqrt{\gamma(\|u_{[0,T]}\|)}\right)\nonumber\\
&\geq\lambda\sqrt{\frac{\underline a}{2\overline a}}\left\| x_{0} \right\| +\frac{1}{2}\sqrt{\frac{\underline a}{2\overline a}}\mu(\|u_{[0,T]}\|), \label{eq-31}
\end{align}
where we used \eqref{eq-39bis} to obtain the last bound. On the other hand, inequality \eqref{GrindEQ__4_} implies that, for almost all $t\ge 0$,
\begin{equation*} 
\frac{d}{dt} \left(x(t)^\top Px(t)\right)\ge -2\sigma\left(\left\| x_{t} \right\| ^{2} +\gamma \left(\left\| u_{[0,t]}\right\| \right) \right).
\end{equation*} 
Recalling that $t_i\leq i(\Delta+\varepsilon)+\Delta$, it holds from the definition of $T$ in \eqref{eq-34} that $t_{i} +\varepsilon \le T$ for all $i\in\{0,\ldots,q-1\}$. It follows from \eqref{GrindEQ__21_} that, for almost all $t\in \left[t_{i} ,t_{i} +\varepsilon \right]$,
\begin{align*} 
\frac{d}{d\, t} \left(x(t)^\top Px(t)\right)&\ge -\frac{4\sigma\overline a}{\underline a}\|x_0\|^2-2\sigma\left(1+\frac{2 T}{\underline a}\right)\gamma(\|u_{[0,t]}\|)\\
&\ge -\frac{4\sigma\overline a}{\underline a}\|x_0\|^2-2\sigma\left(1+\frac{2 T}{\underline a}\right)\gamma(\|u_{[0,T]}\|).
\end{align*} 
By integration, we get for all $i\in\{0,\ldots,q-1\}$ and all $t\in \left[t_{i} ,t_{i} +\varepsilon \right]$ that
\begin{equation*} 
x(t)^\top Px(t)\ge x(t_{i} )^\top Px(t_{i} )-\frac{4\varepsilon\sigma\overline a}{\underline a}\|x_0\|^2-2\sigma\varepsilon\left(1+\frac{2 T}{\underline a}\right)\gamma(\|u_{[0,T]}\|).
\end{equation*} 
It follows that, for all $i\in\{0,\ldots,q-1\}$ and all $t\in \left[t_{i} ,t_{i} +\varepsilon \right]$,
\begin{equation*} 
\left|x(t)\right|^{2}  \ge \frac{p_m  }{p_M  } \left|x(t_{i} )\right|^{2} -\frac{4\varepsilon\sigma\overline a}{p_M\underline a}\|x_0\|^2-\frac{2\sigma\varepsilon}{p_M}\left(1+\frac{2 T}{\underline a}\right)\gamma(\|u_{[0,T]}\|).
\end{equation*} 
Combining this with \eqref{eq-31} yields 
\begin{align*} 
\left|x(t)\right|^{2}  \ge&\, \frac{p_m  }{p_M  } \left(\lambda \sqrt{\frac{\underline a}{2\overline a}}\left\| x_{0} \right\| +\frac{1}{2}\sqrt{\frac{\underline a}{2\overline a}}\mu(\|u_{[0,T]}\|)\right)^2 \\
&\,-\frac{4\varepsilon\sigma\overline a}{p_M\underline a}\|x_0\|^2-\frac{2\sigma\varepsilon}{p_M}\left(1+\frac{2 T}{\underline a}\right)\gamma(\|u_{[0,T]}\|)\\
\geq &\,\frac{p_m \underline a }{p_M  \overline a} \left( \frac{\lambda^2}{2}\|x_0\|^2+\frac{1}{8}\mu(\|u_{[0,T]}\|)^2\right)\\
&\,-\frac{4\varepsilon\sigma\overline a}{p_M\underline a}\|x_0\|^2-\frac{2\sigma\varepsilon}{p_M}\left(1+\frac{2 T}{\underline a}\right)\gamma(\|u_{[0,T]}\|)\\
\geq &\, \frac{1}{p_M}\left( \frac{p_m\lambda^2\underline a}{2\overline a}-\frac{4\varepsilon\sigma\overline a}{\underline a}\right)\|x_0\|^2 \\
&\,+\frac{1}{p_M}\left( \frac{ p_m\underline a}{8\overline a}\mu(\|u_{[0,T]}\|)^2-2\sigma\varepsilon\left(1+\frac{2 T}{\underline a}\right)\gamma(\|u_{[0,T]}\|)\right),
\end{align*} 
for all $t\in[0,T]$. Noticing that the intervals $\left(t_{i} ,t_{i} +\varepsilon \right)$, $i\in\{0,...,q-1\}$, are disjoint subsets of $[0,T]$, we conclude that
\begin{align}
\int _{0}^{T}\left|x(\tau)\right|^{2} d\tau \ge&\, \sum _{i=0}^{q-1}\int _{t_{i} }^{t_{i} +\varepsilon }\left|x(\tau)\right|^{2} d\tau \nonumber\\
\geq &\, \frac{q\varepsilon}{p_M}\left( \frac{p_m\lambda^2\underline a}{2\overline a}-\frac{4\varepsilon\sigma\overline a}{\underline a}\right)\|x_0\|^2 \label{eq-49}\\
&\,+\frac{q\varepsilon}{p_M}\left( \frac{p_m\underline a }{8\overline a}\mu(\|u_{[0,T]}\|)^2-2\sigma\varepsilon\left(1+\frac{2 T}{\underline a}\right)\gamma(\|u_{[0,T]}\|)\right). \nonumber
\end{align}
Observe that $x_0\neq 0$. Indeed, if $x_0$ was zero, then we would have from \eqref{GrindEQ__23_}  and \eqref{eq-39bis} that $\|x_T\|\le \sqrt{\gamma(\|u_{[0,T]}\|)}\sqrt{2T/\underline a}\leq \mu(\|u_{[0,T]}\|)/2$, which is incompatible with \eqref{eq-48}. Based on this, we finally get from \eqref{eq-33}, \eqref{eq-39} and \eqref{eq-49} that
\begin{align*}
\int _{0}^{T}\left|x(\tau)\right|^{2} d\tau > \frac{2\overline a}{\underline a}\|x_0\|^2 +\frac{2T}{a}\gamma(\|u_{[0,T]}\|).
\end{align*}
Confronting this lower bound with the upper bound provided in \eqref{GrindEQ__26_}, we obtain a contradiction. Thus, given any $x_{0} \in \mathcal W^n$ and any $u\in \mathcal U^m$, estimate \eqref{eq-***} holds with the constant $T$ given in \eqref{eq-34}, the function $\mu\in\mathcal N$ proposed in \eqref{eq-38}, and any $\lambda\in(0,1)$ satisfying \eqref{eq-33} for the constants introduced in \eqref{eq-34}. Exp-ISS then directly follows from Lemma \ref{lem-expISS}, by invoking \eqref{GrindEQ__23_} and \eqref{eq-***}. Finally, in view of \eqref{eq-38}, the obtained gain $\mu$ is indeed linear if $\gamma$ is a square function.


\section{Conclusion}

We have presented two growth conditions under which a point-wise dissipation is enough to establish exp-ISS for time-delay systems. These results are new even in the input-free case and significantly extend the class of systems that can be addressed by a point-wise dissipation. They also show that a positive quadratic term involving the whole state history function can be added to the LKF derivative without compromising exp-ISS. Although probably conservative, explicit estimates are provided on the strength of this additional term. We have shown with an example how this extra term can be useful to conduct robustness analysis with respect to modeling uncertainties.

Despite their apparent similarities, the proofs of Theorems \ref{theo-expISS2} and \ref{theo-expISS-left} are radically different: the former relies on the explicit construction of a coercive LKF with history-wise dissipation, whereas the latter relies on a solutions-based contradiction argument. We have not been able to unify the two proofs which, we believe, could provide an interesting direction of research in order to gain insights on these conditions and envision possible ways to further improve them.

Another obvious line of research would be to investigate whether exp-ISS (or GES) can be established under a point-wise dissipation without any growth restrictions. This question remains open and is in line with the conjecture proposed in \cite{CHPEMACH17} that ISS can be established through a point-wise dissipation.

\section*{References}
\bibliography{refs}

\end{document}